\documentclass{article}
\usepackage{amsmath,amsthm,amsopn,amstext,amscd,amsfonts,amssymb,verbatim}
\usepackage{graphicx}
\usepackage{natbib}

\def\cov{\mathop{\rm Cov}\nolimits}

\def\var{\mathop{\rm Var}\nolimits}

\renewenvironment{abstract}
                 {\vspace{6pt}
                  \begin{center}
                  \begin{minipage}{5in}
                  \centerline{\textbf{Abstract}}
                  \noindent\ignorespaces
                 }
                 {\end{minipage}\end{center}}
\newtheorem{thm}{\textbf{Theorem}}[section]

\theoremstyle{definition}

\newtheorem{rem}{\textbf{Remark}}[section]

\title{\Large \textbf{Multivector variate distributions: An application in Finance}}
\author{
  \textbf{\normalsize Jos\'e A. D\'{\i}az-Garc\'{\i}a} \thanks{Corresponding author\newline
   {\bf Key words.}  Bimatrix variate, multivector variate, matrix variate, random vector,
   multivariate elliptical distributions, Kotz distribution.\newline
    2000 Mathematical Subject Classification. 62E15; 60E05}\\
  {\normalsize Universidad Aut\'onoma de Chihuahua, Facultad de Zootecnia y Ecolog\'{\i}a,}\\
  {\normalsize Perif\'erico Francisco R. Almada Km 1, Zootecnia,  33820 Chihuahua, Chihuahua, M\'exico,}\\
  {\normalsize  E-mail: \texttt{jadiaz@uach.mx}}\\
   \textbf{\normalsize Francisco J. Caro-Lopera}\\
  {\normalsize Departament of Basic Sciences, Universidad de Medell\'{\i}n, Medell\'{\i}n, Colombia,} \\
  {\normalsize E-mail: \texttt{fjcaro@udem.edu.co}} \\
  \textbf{\normalsize Fredy O. P\'erez Ram\'{\i}rez}\\
  {\normalsize Faculty of Engineering, Universidad de Medell\'{\i}n, Carrera 87 No.30-65, Medell\'{\i}n, Colombia} \\
  {\normalsize E-mail: \texttt{foperez@udem.edu.co}}
}
\date{}
\begin{document}
\maketitle
\begin{abstract}
A new family of multivariate distributions, which shall be termed multivector variate distributions,
based in the family of the multivariate contoured elliptically distribution is proposed. Several
particular cases of multivector variate distributions are obtained and a number of published multivariate
distributions in another contexts are found as simple corollaries. An application of interest in finance
is full derived and compared with the traditional methods.
\end{abstract}

\section{\large Introduction}\label{sec:1}

Models of phenomena in nature are intrinsically complex and involve a priori several random variables
which depend each other. The researcher expect to discover or model such laws in a parametric or non
parametric way. With observational data the expert tries to estimates the model in order to fit or
predict different stages of interest, however the mathematics behind a model allowing the stochastic
dependency of the variables becomes complex very fast. Then some ideal assumptions must be made involving
independent variables. To clarify this point, we consider a simple example. Consider to model
probabilistically the precipitation and the streamflow, based on the knowledge of a probabilistic
distribution for each variate. The last fact is usually simple to determine, because the expert can
explain very well the marginal distributions of the models, however, the joint distribution of those two
variables is difficult to handle. Then the behavior of both variables in a random bivariate is leaved to
the artificial explanation under independency.  In our example, it is clear that any  occurrence of any
event for streamflow depends on the occurrence of precipitation.

The literature has tackle the problem by constructing a dependent joint distribution via some
transformation, but it is based on independent variables, see for example \citet{n:07} and \citet{n:13}.
That approach based in the knowledge of the marginal distributions is very interesting for application;
in fact the above mention example in hydrology and similar ones in finance have motivated mathematical
developments in such area, see also \citet{cn:84} and \citet{ol:03}, unfortunately, all of them are based
in some way in independent variables.

From the theoretical point of view, the referred multivariate case has been generalizsed to the matrix
case, by the so called bimatrix variate distributions, see \citet{brea:11} and
\citet{dggj:10a,dggj:10b,dggj:11}. The reader is also referred to \citet{e:11}, which provides a complete
revision of the works about such complex problem.

We  quote again that the best advantage of the model consists of knowing the marginal distribution of
each variable, a simple requirement to fulfill in applied sciences, given the records of univariate
random research.  The main problem to solve consists of proposing a probabilistic dependent model with
such marginal distributions. We set the solution as follows: Section \ref{sec:2} gives some results on
Jacobians, integrals and elliptical distributions. Then Section \ref{sec:3} derives the family of joint
distributions of several probabilistic dependent random vectors with known marginal distributions. Such
distribution which is indexed by the general family of elliptical contoured distributions will be termed
''multi-vector variate distribution". As corollaries, several particular distributions of the literature
are derived straightforwardly, and a number of new distributions and versatile mixtures of distributions
are provided in  Section \ref{sec:4}. In Section \ref{sec:5}, distributions studied in  previous section
are extended to more general values of their parameters. Finally an applications in finance is full
detailed in Section \ref{sec:6}.

\section{\large Preliminary results}\label{sec:2}

In this section we shown some results about Jacobians, integration and multivariate elliptical contoured
distributions necessary for the presentation of the results of this work. A detailed study of these and
other results can be consulted in \citet{fz:90}, \citet{fkn:90} and \citet{mh:05}.

First we propose some notations: if $\mathbf{A}\in \Re^{n \times m}$ denotes a \emph{matrix}, this is,
$\mathbf{A}$ have $n$ rows and $m$ columns, then $\mathbf{A}'\in \Re^{m \times n}$ denotes its
\emph{transpose matrix}, and if $\mathbf{A}\in \Re^{n \times n}$ has an \emph{inverse}, it shall be
denoted by $\mathbf{A}^{-1} \in \Re^{n \times n}$. An \emph{identity matrix} shall be denoted by
$\mathbf{I}\in \Re^{n \times n}$, to specified the size of the identity, we will use $\mathbf{I}_{n}$. A
\emph{null matrix} shall be denoted as $\mathbf{0} \equiv \mathbf{0}_{n \times m} \in \Re^{n\times m}$.
$\mathbf{A}\in \Re^{n \times n}$ is a \emph{symmetric matrix} if $\mathbf{A} = \mathbf{A}'$ and if all
their eigenvalues are positive then $\mathbf{A}$ is \emph{positive definite matrix},  which shall be
denoted as $\mathbf{A} > \mathbf{0}$. The vectors shall be denoted also in bold style but in lower case.
The Euclidean norm of a vector $\mathbf{x}=(x_{1},\ldots,x_{n})'\in\Re^{n}$ will be denoted as
$||\mathbf{x}||$ and is given by $\sqrt{\mathbf{x}'\mathbf{x}}=\sqrt{\sum_{i=1}^{n}x_{i}^{2}}$.

\bigskip

From \citet[Theorem 2.1.3, p.55]{mh:05}, we have,
\begin{thm}\label{teo1}
For the following transformation from rectangular coordinates $\mathbf{x}' = (x_{1}, \dots, x_{n})$ to
polar coordinates $(s, \mathbf{h})' = (s, \theta_{1}, \dots, \theta_{n-1})$:
\begin{eqnarray*}
  x_{1} &=& s \sin \theta_{1} \sin \theta_{2}\cdots\sin \theta_{n-2} \sin \theta_{n-1}\\
  x_{2} &=& s \sin \theta_{1} \sin \theta_{2}\cdots\sin \theta_{n-2} \cos \theta_{n-1}\\
  x_{3} &=& s \sin \theta_{1} \sin \theta_{2}\cdots\sin \theta_{n-3} \cos \theta_{n-2}\\
  \vdots &=& \vdots \\
  x_{n-1} &=& s \sin \theta_{1}\cos \theta_{2} \\
  x_{n} &=& s \cos \theta_{1},
\end{eqnarray*}
where, $s > 0$, $\theta_{i} \in (0, \pi]$, $i = 1, \dots, n-2$ and $\theta_{n-1} \in (0, 2\pi]$. We have
$$
 \bigwedge_{i = 1}^{n}dx_{i} = s^{n-1}\prod_{i=1}^{n-2} \sin^{n-i-1} \theta_{i}
 \left(\bigwedge_{i=}^{n-1}d\theta_{i}\right)\wedge ds.
$$
\end{thm}
In this results is defined $s^{2} = \|\mathbf{x}\|^{2}$. Then if we define $s = r^{2}$, we obtain
$$
 \bigwedge_{i = 1}^{n}dx_{i} = 2^{-1} r^{n/2-1}\prod_{i=1}^{n-2} \sin^{n-i-1} \theta_{i}
 \left(\bigwedge_{i=}^{n-1}d\theta_{i}\right)\wedge dr.
$$
Denoting
$$
 \prod_{i=1}^{n-2} \sin^{n-i-1} \theta_{i} \left(\bigwedge_{i=}^{n-1}d\theta_{i}\right),
$$
as $(\mathbf{h}'d\mathbf{h})$, we obtain that
$$
 \bigwedge_{i = 1}^{n}dx_{i} = 2^{-1}r^{n/2-1}(\mathbf{h}'d\mathbf{h})\wedge dr.
$$
Moreover, from the proof of \citet[Theorem 1.5.5, pp. 36-37 and Thoerem 2.1.15, p. 70]{mh:05}, we have
$$
  \int_{\mathbf{h} \in \mathcal{V}_{1,n}} (\mathbf{h}'d\mathbf{h}) = \frac{2 \pi^{n/2}}{\Gamma[n/2]},
$$
where in general, $\mathcal{V}_{m,n}$ denotes the \emph{Stiefel manifold}, the space of all matrices
$\mathbf{H}_{1} \in \Re^{n \times m}$ ($n \geq m$) with orthogonal columns, so that
$\mathbf{H}'_{1}\mathbf{H}_{1} = \mathbf{I}_{m}$. Thus
$$
  \mathcal{V}_{m,n} = \{\mathbf{H}_{1} \in \Re^{n \times m}| \mathbf{H}'_{1}\mathbf{H}_{1} =
  \mathbf{I}_{m}\}.
$$
Note that, when $m = 1$, then
$$
  \mathcal{V}_{1,n} = \{\mathbf{h} \in \Re^{n}: \|\mathbf{h}\|^{2} = 1\},
$$
the unit sphere in $\Re^{n}$. This is, of course, an $(n-1)$-dimensional surface in $\Re^{n}$, see
\citet[Section 2.1.4, pp. 67-72]{mh:05}.\bigskip

Thus we can propose an important Jacobian.

\begin{thm}
Let $\mathbf{y}=\left(1-||\mathbf{x}||^{2}\right)^{-\frac{1}{2}}\mathbf{x}$, with $\mathbf{x},
\mathbf{y}\in \Re^{n}$, then
\begin{equation*}
    (d\mathbf{y})=\left(1-||\mathbf{x}||^{2}\right)^{-\left(\frac{n}{2}+1\right)}(d\mathbf{x}).
\end{equation*}
\end{thm}
\proof Define $v=||\mathbf{y}||^{2}=\left(1-||\mathbf{x}||^{2}\right)^{-1}||\mathbf{x}||^{2}$ and
$w=||\mathbf{x}||^{2}$, then by \citet[Th. 2.1.14]{mh:05},
\begin{equation*}
    (d\mathbf{y})=2^{-1}v^{\frac{n}{2}-1}(dv)\wedge(\mathbf{h}'d\mathbf{h}),  \quad \mathbf{h}\in \mathcal{V}_{1,n}
\end{equation*}
and
\begin{equation*}
    (d\mathbf{x})=2^{-1}w^{\frac{n}{2}-1}(dw)\wedge(\mathbf{g}'d\mathbf{g}), \quad \mathbf{g}\in \mathcal{V}_{1,n},
\end{equation*}
then
\begin{equation}\label{eq1}
    (dv)=2v^{-\left(\frac{n}{2}-1\right)}(d\mathbf{y})\wedge(\mathbf{h}'d\mathbf{h}), \quad \mathbf{h}\in \mathcal{V}_{1,n}
\end{equation}
and
\begin{equation}\label{eq2}
    (dw)=2w^{-\left(\frac{n}{2}-1\right)}(d\mathbf{x})\wedge(\mathbf{g}'d\mathbf{g}), \quad \mathbf{g}\in \mathcal{V}_{1,n},
\end{equation}
also $v=(1-w)^{-1}w=
(1-w)^{-1}(1-(1-w))=(1-w)^{-1}-1$, thus
\begin{equation}\label{eq3}
  (dv)=\mbox{abs}\left(\frac{(dw)}{(1-w)^{2}}\right).
\end{equation}
Uniqueness of $(\mathbf{h}'d\mathbf{h})$  and substitution of (\ref{eq1}) and (\ref{eq2}) in (\ref{eq3})
imply that $(d\mathbf{y})=\frac{v^{\frac{n}{2}-1}}{(1-w)^{2}}(d\mathbf{x})$. Thus  using
$\mathbf{y}=\left(1-||\mathbf{x}||^{2}\right)^{-\frac{1}{2}}\mathbf{x}$, $v=||\mathbf{y}||^{2}$ and
$w=||\mathbf{x}||^{2}$, the required result is obtained. \qed

Next we provide from \citet{mh:05} a summary about elliptically  contoured distribution theory; see also
\citet{fz:90} and \citet{fkn:90}.  Let $\mathbf{y}=(y_{1},\ldots,y_{n})' \in \Re^{n}$, then $\mathbf{y}$
is said to have a \emph{multivariate elliptically contoured distribution} if its characteristic function
takes the form
$$
  \phi_{\mathbf{y}}(\mathbf{t})=\exp(i\mathbf{t}'\boldsymbol{\mu})\psi(\mathbf{t}'\boldsymbol{\Sigma}\mathbf{t}),
$$
where
 $\mathbf{t}$ and $\boldsymbol{\mu}$ are n-dimensional vectors and $\psi:[0,\infty)\rightarrow \Re$,
$|\psi(t)|\leq 1$ for $t\in\Re_{0}^{+}$. It implies that all marginal distributions are elliptical and
all marginal density functions of dimension $k<n$ have the same functional form. For example, take
$\mathbf{t}=(\mathbf{t}_{1}'|\mathbf{0}')'$, $\mathbf{y}=\left(
                   \begin{array}{c}
                     \mathbf{y}_{1} \\
                     \mathbf{y}_{2} \\
                   \end{array}
                 \right)
$, $\boldsymbol{\mu}=\left(
                   \begin{array}{c}
                     \boldsymbol{\mu}_{1} \\
                     \boldsymbol{\mu}_{2} \\
                   \end{array}
                 \right)
$, $\boldsymbol{\Sigma}=\left(
      \begin{array}{cc}
        \boldsymbol{\Sigma}_{11} & \boldsymbol{\Sigma}_{12} \\
        \boldsymbol{\Sigma}_{21} & \boldsymbol{\Sigma}_{22} \\
      \end{array}
    \right)
$, where $\mathbf{t}_{1}$, $\mathbf{y}_{1}$ and $\boldsymbol{\mu}_{1}$ are $k\times 1$ and
$\boldsymbol{\Sigma}_{11}$ is $k\times k$, then via the characteristic function we have that
$\mathbf{y}_{1}\sim \mathcal{E}_{k}(\boldsymbol{\mu}_{1},\boldsymbol{\Sigma}_{11},\psi)$. Provided they
exist, $E(\mathbf{y})=\boldsymbol{\mu}$ and $\cov(\mathbf{y})=c\boldsymbol{\Sigma}$, where
$c=-2\psi'(0)$. Thus, all distributions in the class
$\mathcal{E}_{n}(\boldsymbol{\mu},\boldsymbol{\Sigma},\psi)$ have the same mean $\boldsymbol{\mu}$ and
the same correlation matrix $\mathbf{P}=(\rho_{ij})$, with
$$
  \rho_{ij}=\frac{\cov(y_{i},y_{j})}{\sqrt{\var(y_{i})}\sqrt{\var(y_{j})}}=\frac{c\sigma_{ij}}{\sqrt{c\sigma_{ii}}
  \sqrt{c\sigma_{jj}}}=\frac{\sigma_{ij}} {\sqrt{\sigma_{ii}}\sqrt{\sigma_{jj}}}.
$$
Note that if any marginal distribution is normal, then $\mathbf{y}$ is normal. Observe also that if
$\mathbf{y}=(y_{1},\ldots,y_{n})'\sim \mathcal{N}_{n}(\boldsymbol{\mu},\boldsymbol{\Sigma})$, and
$\boldsymbol{\Sigma}$ is diagonal, then $y_{1},\ldots,y_{n}$ are all independent. In general, let
$\mathbf{y}=(y_{1},\ldots,y_{n})'\sim \mathcal{E}_{n}(\boldsymbol{\mu},\boldsymbol{\Sigma},\psi)$ and
$\boldsymbol{\Sigma}$ is diagonal, if $y_{1},\ldots,y_{n}$ are all independent, then $\mathbf{y}$ is
normal.

Now, the density function of the elliptical contoured n-dimensional vector $\mathbf{y}$, with respect to
the Lebesgue measure in $\Re^{n}$, is given by
\begin{equation}\label{e}
    dF_{\mathbf{y}}(\mathbf{y})=\frac{1}{|\boldsymbol{\Sigma}|^{1/2}}
    h\left\{\left[(\mathbf{y}-\boldsymbol{\mu})'\boldsymbol{\Sigma}^{-1}(\mathbf{y}-
    \boldsymbol{\mu})\right]\right\} (d\mathbf{y}),
\end{equation}
where  $\boldsymbol{\mu} \in \Re^{n}$, $\boldsymbol{\Sigma} \in \Re^{n\times n}$,
$\boldsymbol{\Sigma}>\mathbf{0}$  and $(d\mathbf{Y})=\bigwedge_{i=1}^{n}dy_{i}, i=1,\ldots,n$. The
function $h: \Re \rightarrow [0,\infty)$ is termed the generator function and satisfies $\int_{0}^\infty
u^{n-1}h(u^2)du < \infty$. Such a distribution is denoted by $\mathbf{y}\sim
\mathcal{E}_{n}(\boldsymbol{\mu},\boldsymbol{\Sigma}, h)$.

This class of distributions includes the Normal, t-Student (which includes the t and Cauchy
distributions), Laplace, Bessel and Kotz, among other distributions. Some specific matrix-variate
Elliptical densities are presented in Table 1, also see \citet{fz:90}.

\def\baselinestretch{1}
\begin{table}[!ht]\label{table0}
\caption{Explicit forms of the densities of the elliptical
distributions.}
\medskip
\hspace{1.2cm}
\begin{minipage}[t]{360pt}
\begin{small}
\begin{center}
\begin{tabular}{|| l | c||}
\hline \hline Elliptical law &  Density\footnote{\scriptsize Where
$W =  (\mathbf{y} -
\boldsymbol{\mu})'\boldsymbol{\Sigma}^{-1}(\mathbf{y} - \boldsymbol{\mu})$ and $\varrho = |\boldsymbol{\Sigma}|^{1/2}$} \\
\hline \hline&\\
Pearson VII\footnote{\scriptsize Particularly, when $q = (n+r)/2$
we obtain the matrix-variate $t$-distribution. And when $r = 1$ the
matrix-variate Cauchy distribution is obtained.}
& $\displaystyle\frac{\Gamma[q]}{(r \pi)^{n/2} \Gamma[q-n/2] \varrho}\left ( 1 + \displaystyle\frac{W}{r}\right)^{-q}$\\
&{\footnotesize $r>0$ and $q > n/2$}\\
\hline  &\\
Kotz type\footnote{\scriptsize If we take $q=s=1$ and $r=1/2$ we
obtain the matrix-variate normal distribution.} &
$\displaystyle\frac{s r^{(2q+n-2)/2s} \ \Gamma[n/2]}{\pi^{n/2}
\Gamma\left[(2q+n-2)/2s \right] \varrho} \ W^{q-1}\exp\left(-r W^{s}\right)$\\
&{\footnotesize$r, s > 0$ and $2q + n > 2$}\\
\hline  &\\
Bessel\footnote{\scriptsize With
$$
    K_{q}(z) = \displaystyle\frac{\pi}{2} \displaystyle\frac{I_{-q}(z)-I_{q}(z)}{\sin(q\pi)}, \ |\arg(z)| < \pi,
$$
with $q$ an integer, is the modified Bessel function of the third
kind and
$$
    I_{q}(z) = \displaystyle\sum^{\infty}_{k=0} \displaystyle\frac{1}{k! \Gamma[k+q+1]}
\left(\displaystyle\frac{z}{2}\right)^{q+2k}, \ |z| < \infty, \
|\arg(z)| < \pi.
$$
Also, note that for $q = 0$ and $r = \sigma/\sqrt{2}$, $\sigma>0$,
the matrix-variate
Laplace distribution is obtained.}%
& $\displaystyle\frac{W^{1/2}}{2^{q+n-1} \ \pi^{n/2} \ r^{n+q} \
\Gamma[q + n/2] \ \varrho}
    \ K_{q}\left(\displaystyle\frac{W^{1/2}}{r}\right)$ \\
&{\footnotesize$r>0$ and $q > -n/2$}\\
\hline  &\\
Pearson II& $\displaystyle\frac{\Gamma[q+1+n/2]}{\pi^{n/2} \
\Gamma[q+1]\ \varrho} (1 -
W)^{q}$\\
&{\footnotesize$q \in \Re$ and $W \leq 1$}\\
\hline \hline
\end{tabular}
\end{center}
\end{small}
\end{minipage}
\end{table}

When $\boldsymbol{\mu}=\mathbf{0}_{n}$, $\boldsymbol{\Sigma}= \mathbf{I}_{n}$ , such distribution is
termed \emph{vector variate spherical distribution} and shall be denoted as $\mathbf{Y} \sim
\mathcal{E}_{n}(\mathbf{0}, \mathbf{I}_{n}, h)$. In other words, $\mathbf{y}$ is said to have a spherical
distribution if $\mathbf{y}$ and $\mathbf{H}\mathbf{y}$ have the same distribution for all $n\times n$
orthogonal matrices $\mathbf{H}$. Moreover, the spherical density depends on $\mathbf{y}$ only through
the value of $\mathbf{y}'\mathbf{y}$.

\bigskip

Now, assume that $\mathbf{x}\sim \mathcal{E}_{n}(\mathbf{0},\mathbf{I},h)$, i.e.
\begin{equation*}
dF_{\mathbf{x}}(\mathbf{x})=h(||\mathbf{x}||^{2})(d\mathbf{x}).
\end{equation*}
Consider the change of rectangular coordinates into polar coordinates, $\mathbf{x}\rightarrow
(s,\mathbf{g})$, then by Theorem \ref{teo1}
\begin{equation*}
    dF_{s,\mathbf{g}}(s,\mathbf{g})=2^{-1}s^{n/2-1}h(s)(ds)\wedge(\mathbf{g}'d\mathbf{g}), \quad s>0, ||\mathbf{g}||^{2}=1.
\end{equation*}
Integration over $\mathbf{g}\in \mathcal{V}_{1,n}$ gives
\begin{equation*}
    dF_{s}(s)=\frac{\pi^{n/2}}{\Gamma[n/2]}s^{n/2-1}h(s)(ds),s>0.
\end{equation*}
Set $r=s^{1/2}$, then $(ds)=2r(dr)$ and
\begin{equation*}
    dF_{r}(r)=\frac{2\pi^{n/2}}{\Gamma[n/2]}s^{n-1}h(r^{2})(dr).
\end{equation*}
Observe also that
\begin{equation*}
    1=\int_{\mathbf{x}\in \Re^{n}}h(||\mathbf{x}||^{2})(dx),
\end{equation*}
so by the change of variable $z=as$, with $(ds)=\frac{(dz)}{a}$, we have that
\begin{equation*}
    1=\int_{s>0}\frac{\pi^{n/2}}{\Gamma[n/2]}s^{n/2-1}h(s)(ds)=
    \frac{\pi^{n/2}}{\Gamma[n/2]}\int_{z>0}
    \left(\frac{z}{a}\right)^{n/2-1}h\left(\frac{z}{a}\right)\frac{(dz)}{a}.
\end{equation*}
Thus we get that
\begin{equation}\label{eqintgen}
    \int_{z>0}
    z^{n/2-1}h\left(\frac{z}{a}\right)(dz)=\frac{a^{n/2}\Gamma[n/2]}{\pi^{n/2}}.
\end{equation}

\bigskip

\section{\large Multivector variate theory}\label{sec:3}

In this apart we derive a number of multivector variate distribution. The technique can be replicated
with several mixtures or marginal distributions, which depends on the particular phenomena of interest
and a previous knowledge of the marginal distributions.

\subsection{Multivector variate elliptical distribution}
Assume that $\mathbf{x}\sim \mathcal{E}_{n}(\boldsymbol{\mu},\boldsymbol{\Sigma},h)$, with
$\boldsymbol{\Sigma}>\mathbf{0}$, then
\begin{equation*}
    dF_{\mathbf{x}}(\mathbf{x})=|\boldsymbol{\Sigma}|^{-\frac{1}{2}}h\left((\mathbf{x}-\boldsymbol{\mu})'\boldsymbol{\Sigma}'
    (\mathbf{x}-\boldsymbol{\mu})\right)(d\mathbf{x}).
\end{equation*}
Now consider the following partition of $\mathbf{x}$, $\boldsymbol{\mu}$ and $\boldsymbol{\Sigma}$
\begin{eqnarray*}
  &&\mathbf{x} =\left(\begin{array}{c}
                    \mathbf{x}_{1} \\
                    \vdots \\
                    \mathbf{x}_{k}
                  \end{array}\right) \in\Re^{n}, \boldsymbol{\mu} = \left(\begin{array}{c}
                    \boldsymbol{\mu_{1}} \\
                    \vdots \\
                    \boldsymbol{\mu_{k}}
                  \end{array}\right), \mathbf{x}_{i}, \boldsymbol{\mu}_{i}\in\Re^{n_{i}},\sum_{i=1}^{k}n_{i}=n,
   \\
   &&\boldsymbol{\Sigma}=\left(
                           \begin{array}{cccc}
                             \boldsymbol{\Sigma}_{11} & \mathbf{0} & \cdots & \mathbf{0} \\
                             \mathbf{0} & \boldsymbol{\Sigma}_{22} & \cdots & \mathbf{0} \\
                             \vdots & \vdots & \ddots & \vdots \\
                             \mathbf{0} & \mathbf{0} & \cdots & \boldsymbol{\Sigma}_{kk} \\
                           \end{array}
                         \right), \boldsymbol{\Sigma}_{ii}>\mathbf{0}, \boldsymbol{\Sigma}_{ii}\in\Re^{n_{i}\times n_{i}}, i=1,\ldots,k.
\end{eqnarray*}
Then we have
\begin{equation}\label{eqmultielliptical}
    dF_{\mathbf{x}_{1},\ldots,\mathbf{x}_{k}}(\mathbf{x}_{1},\ldots,\mathbf{x}_{k})=
    \prod_{i=1}^{k}|\boldsymbol{\Sigma}_{ii}|^{-1/2}h\left(\sum_{i=1}^{k}
    (\mathbf{x}_{i}-\boldsymbol{\mu}_{i})'\boldsymbol{\Sigma}_{ii}^{-1}(\mathbf{x}_{i}-\boldsymbol{\mu}_{i})\right)
    \bigwedge_{i=1}^{k}(d\mathbf{x}_{i}).
\end{equation}
This joint distribution define a multivector variate which shall be termed \emph{multivector variate
elliptical distribution}. It will be denoted as
\begin{equation*}
    \mathbf{x} =\left(\begin{array}{c}
                    \mathbf{x}_{1} \\
                    \vdots \\
                    \mathbf{x}_{k}
                  \end{array}
                  \right) \sim \mathcal{ME}_{n}(n_{1},\ldots,n_{k};\boldsymbol{\mu_{1}},\ldots,\boldsymbol{\mu_{k}};
                  \boldsymbol{\Sigma_{11}}, \ldots,\boldsymbol{\Sigma_{kk}};h).
\end{equation*}

\subsection{Multivector variate log-elliptical distribution}

Now assume that $\mathbf{x}=\log \mathbf{v}$, $\mathbf{x}, \mathbf{v}\in\Re^{n}$, and $\log
\mathbf{v}=(\log v_{1},\ldots,\log v_{n})'$. If $\mathbf{x}\sim
\mathcal{E}_{n}(\boldsymbol{\mu},\boldsymbol{\Sigma},h)$, then  $\mathbf{v}$ is said to have a
log-elliptical distribution, see \citet{fkn:90}.

Consider that $\mathbf{x} =\log \mathbf{v}\sim \mathcal{ME}_{n}(n_{1},\ldots,n_{k};
\boldsymbol{\mu_{1}},\ldots,\boldsymbol{\mu_{k}}; \boldsymbol{\Sigma_{11}},
\ldots,\boldsymbol{\Sigma_{kk}};h)$, where $\mathbf{x}_{i}=\log \mathbf{v}_{i}$, $i=1,\ldots,k$, then we
say that $\mathbf{v}=(\mathbf{v}_{1},\ldots,\mathbf{v}_{k})'$ has a \emph{multivector log-elliptical
distribution}, denoted as
 \begin{equation*}
    \mathbf{v} =\left(\begin{array}{c}
                    \mathbf{v}_{1} \\
                    \vdots \\
                    \mathbf{v}_{k}
                  \end{array}\right)\sim \mathcal{MLE}_{n}(n_{1},\ldots,n_{k};\boldsymbol{\mu_{1}},\ldots,\boldsymbol{\mu_{k}};
                  \boldsymbol{\Sigma_{11}}, \ldots,\boldsymbol{\Sigma_{kk}};h).
\end{equation*}

Moreover, making the change of variable $(\mathbf{x}_{1},\ldots,\mathbf{x}_{k})'= (\log\mathbf{v}_{1},
\ldots, \log\mathbf{v}_{k})'$, we have that
\begin{eqnarray*}
  \bigwedge_{i=1}^{k}(d\mathbf{x}_{i}) &=& \left(\prod_{i=1}^{n_{1}}v_{i}^{-1}\right)\cdots\left(\prod_{i=1}^{n_{k}}v_{i}^{-1}\right)
  \bigwedge_{i=1}^{k}(d\mathbf{v}_{i}) \\
   &=& \prod_{j=1}^{k}\prod_{i=1}^{n_{j}}v_{i}^{-1}\bigwedge_{i=1}^{k}(d\mathbf{v}_{i}).
\end{eqnarray*}
And by (\ref{eqmultielliptical}) we obtain
$$
    dF_{\mathbf{v}_{1},\ldots,\mathbf{v}_{k}}(\mathbf{v}_{1},\ldots,\mathbf{v}_{k})= \prod_{j=1}^{k}\left[|\boldsymbol{\Sigma}_{jj}|^{-n/2}
    \left(\prod_{i=1}^{n_{j}}v_{i}^{-1}\right)\right]\hspace{6cm}
$$
\begin{equation}\label{eqmultilogelliptical}
   \hspace{4cm}\times h\left(\sum_{i=1}^{k}(\log\mathbf{v}_{i}-\boldsymbol{\mu}_{i})' \boldsymbol{\Sigma}_{ii}^{-1}
    (\log\mathbf{v}_{i}-\boldsymbol{\mu}_{i})\right) \bigwedge_{i=1}^{k}(d\mathbf{v}_{i})
\end{equation}
This distribution is multifunctional because we can  find a number of new models. For example, we can
obtain multivector variate mixed distribution, namely, a part of the multivector $\mathbf{v}$ can follow
a multivector variate elliptical distribution and the remaining part indexed by multivector variate
log-elliptical law. This fact will be denoted as
\begin{equation*}
    \left(
      \begin{array}{c}
        \mathbf{x} \\
        \mathbf{v} \\
      \end{array}
    \right)
=\left(\begin{array}{c}
                    \mathbf{x}_{1} \\
                    \vdots \\
                    \mathbf{x}_{k_{1}}\\
                    \mathbf{v}_{1} \\
                    \vdots \\
                    \mathbf{v}_{k_{2}}
                  \end{array}\right)\sim \mathcal{MELE}_{n}(n_{1},\ldots,n_{k};\boldsymbol{\mu_{1}},\ldots,\boldsymbol{\mu_{k}};
                  \boldsymbol{\Sigma_{11}}, \ldots,\boldsymbol{\Sigma_{kk}};h),
\end{equation*}
with $k_{1}+k_{2}=k$, $\mathbf{x}_{i}\in\Re^{n_{i}}$, $\mathbf{v}_{i}\in\Re^{n_{i}}$,
$n=\sum_{i=1}^{k}n_{i}$ and the corresponding joint distribution is given by
$$
  dF_{\mathbf{x}_{1},\ldots,\mathbf{x}_{k_{1}},\mathbf{v}_{1},\ldots,\mathbf{v}_{k_{2}}}
    (\mathbf{x}_{1},\ldots,\mathbf{x}_{k_{1}},\mathbf{v}_{1},\ldots,\mathbf{v}_{k_{2}}) \hspace{6cm}
$$
$$
    = \prod_{i=1}^{k_{1}}|\boldsymbol{\Sigma}_{ii}|^{-n/2} \prod_{j=1}^{k_{2}}\left[|\boldsymbol{\Sigma}_{jj}|^{-n/2}
    \left(\prod_{r=1}^{n_{j}}v_{r}^{-1}\right)\right] \hspace{2cm}
$$
$$
  \times h\left(\sum_{i=1}^{k_{1}}  (\mathbf{x}_{i}-\boldsymbol{\mu}_{i})'\boldsymbol{\Sigma}_{ii}^{-1}(\log\mathbf{x}_{i}
  -\boldsymbol{\mu}_{i}) \right. \hspace{1cm}
$$
$$\hspace{2.5cm}
   + \left.\sum_{j=1}^{k_{2}}(\log\mathbf{v}_{j}-\boldsymbol{\mu}_{j})'\boldsymbol{\Sigma}_{jj}^{-1}
   (\log\mathbf{v}_{j}-\boldsymbol{\mu}_{j})\right) \bigwedge_{i=1}^{k_{1}}(d\mathbf{x}_{i})
   \bigwedge_{j=1}^{k_{2}}(d\mathbf{v}_{j}).
$$

Note that in this distribution, the variables $\mathbf{x}_{i}$, $i = 1,\dots,k_{1}$ are probabilistically
dependent, the same property follows for the variables $\mathbf{v}_{r}$, $r = 1,\dots,k_{2}$. Moreover,
for all $i = 1,\dots,k_{1}$ and $r = 1,\dots,k_{2}$, the variates $\mathbf{x}_{i}$ and $\mathbf{v}_{r}$
are also probabilistically dependent. This family of distribution can be termed as  \emph{multivector
variate elliptical-log-elliptical distribution}. In particular,  for $k_{1}=k_{2}=1$ this distribution
can be proposed as the joint distribution for modeling two random variables $x$ and $v$ with certain
marginal distribution; for example, $x$ could follow a  $t$-Student distribution and  $v$ a $\log$
$t$-Student distribution; we highlight again that the addressed joint distribution considers a
probabilistic dependency between the two random variables $x$ and $v$.

\subsection{Multivector variate T-distribution (or Pearson type VII)}

Recall that if $\mathbf{x} \sim \mathcal{N}_{n}(\mathbf{0},\mathbf{I})$ and $s^{2}\sim \chi_{\nu}^{2}$
are independent, then $\mathbf{t}=s^{-1}\mathbf{x}$ has a multivariate t distribution.  Now consider that
$$
\left(\begin{array}{c}
                    \mathbf{x}_{0}\\
                    \mathbf{x}_{1} \\
                    \vdots \\
                    \mathbf{x}_{k}
                  \end{array}\right)\sim \mathcal{ME}_{n}(n_{0},\ldots,n_{k};\mathbf{0},\ldots,\mathbf{0};
                  \sigma_{0}^{2}\mathbf{I}_{n_{0}}, \ldots,\sigma_{k}^{2}\mathbf{I}_{n_{k}};h),
$$
and define $s_{0}=||\mathbf{x}_{0}||^{2}$ and $\mathbf{t}_{i}=s_{0}^{-1/2}\mathbf{x}_{i}$,
$i=1,\ldots,k$. Then
\begin{equation*}
   \bigwedge_{i=0}^{k}(d\mathbf{x}_{i})=2^{-1}s_{0}^{n_{0}/2-1}\prod_{i=1}^{k}s_{0}^{n_{i}/2}
   (ds_{0})\wedge(\mathbf{g}'d\mathbf{g})\wedge\left(\bigwedge_{i=1}^{k}(d\mathbf{t}_{i})\right),\quad \mathbf{g}\in \mathcal{V}_{1,n_{0}}.
\end{equation*}
Thus

\begin{equation}\label{eqmultitx}
    dF_{\mathbf{x}_{0},\ldots,\mathbf{x}_{k}} (\mathbf{x}_{0},\ldots,\mathbf{x}_{k})=
    \left(\prod_{i=1}^{k}\sigma_{i}^{-n_{i}}\right) h\left(\sum_{i=0}^{k}\sigma_{i}^{-2}||\mathbf{x}_{i}||^{2}\right)
    \bigwedge_{i=1}^{k}(d\mathbf{x}_{i})
\end{equation}
so
$$
    dF_{s_{0},\mathbf{g},\mathbf{t}_{1},\ldots,\mathbf{t}_{k}}
    (s_{0},\mathbf{g},\mathbf{t}_{1},\ldots,\mathbf{t}_{k})=\frac{s_{0}^{n^{*}/2-1}}
    {2\prod_{i=0}^{k}\sigma_{i}^{n_{i}}} \hspace{6cm}
$$
\begin{equation}\label{eqmultitjoin}
    \hspace{2cm}\times h\left(\sigma_{0}^{-2}s_{0}+s_{0}\sum_{i=1}^{k}\sigma_{i}^{-2}||\mathbf{t}_{i}||^{2}\right)
    (ds_{0})\wedge(\mathbf{g}'d\mathbf{g})\wedge \bigwedge_{i=1}^{k}(d\mathbf{t}_{i}),
\end{equation}
where $n^{*} = n_{0}+n_{1}+\cdots+n_{k}$. Integration over $\mathbf{g}\in\mathcal{V}_{1,n_{0}}$ gives
$$
  dF_{s_{0},\mathbf{t}_{1},\ldots,\mathbf{t}_{k}}(s_{0},\mathbf{t}_{1},\ldots,\mathbf{t}_{k}) =
  \frac{\pi^{n_{0}/2}}{\Gamma\left[n_{0}/2\right]\prod_{i=0}^{k}\sigma_{i}^{n_{i}}} s_{0}^{n^{*}/2-1}
  \hspace{4cm}
$$
\begin{equation}\label{eqmultitjoinst}
  \hspace{2cm} \times h\left(\sigma_{0}^{-2}\left(1+\sigma_{0}^{2}\sum_{i=1}^{k} \sigma_{i}^{-2}
  ||\mathbf{t}_{i}||^{2}\right)s_{0}\right)(ds_{0}) \wedge \left(\bigwedge_{i=1}^{k}(d\mathbf{t}_{i})\right).
\end{equation}

\begin{rem}
According to the definition of this distribution, we shall see that the random variable $s_{0}$ follows a
generalised $\chi^{2}$ distribution (or generalised gamma distribution), see \citet{fz:90}, and
$\mathbf{t} = (\mathbf{t}'_{1}, \dots, \mathbf{t}'_{k})'$ has a multivector variate $t$-distribution (or
Pearson type VII). Also note that  $s_{0}$ and $\mathbf{t}$ are probabilistically dependent. Finally,
this family of distribution defines a \emph{multivector variate $\chi^{2}$ generalisada-$t$
distribution}.
\end{rem}

From (\ref{eqintgen}), the integral with respect to $s_{0}$ is given by
$$
    dF_{\mathbf{t}_{1},\ldots,\mathbf{t}_{k}}
    (\mathbf{t}_{1},\ldots,\mathbf{t}_{k})=\frac{\Gamma\left[n^{*}/2\right]\sigma_{0}^{n^{*}/2}
    }{\Gamma\left[n_{0}/2\right]\pi^{n^{*}/2}\prod_{i=0}^{k}\sigma_{i}^{n_{i}}} \hspace{6cm}
$$
\begin{equation}\label{eqmultitt1}
    \hspace{5cm} \times \left(1+\sigma_{0}^{2}\sum_{i=1}^{k}\sigma_{i}^{-2}||\mathbf{t}_{i}||^{2}\right)^{-n^{*}/2}
    \bigwedge_{i=1}^{k}(d\mathbf{t}_{i}).
\end{equation}
Finally, let $\beta_{i}=\sigma_{i}^{2}/\sigma_{0}^{2}$, then
$$
    dF_{\mathbf{t}_{1},\ldots,\mathbf{t}_{k}}
    (\mathbf{t}_{1},\ldots,\mathbf{t}_{k})=\frac{\Gamma\left[n^{*}/2\right]\prod_{i=1}^{k}\beta_{i}^{-n_{i}/2}}
    {\Gamma\left[n_{0}/2\right]\pi^{n^{*}/2}} \hspace{6cm}
$$
\begin{equation}\label{eqmultitt}
    \hspace{5cm} \times \left(1+\sum_{i=1}^{k}\beta_{i}^{-1}||\mathbf{t}_{i}||^{2}\right)^{-n^{*}/2}
    \bigwedge_{i=1}^{k}(d\mathbf{t}_{i}).
\end{equation}
We refer it as the \emph{multivector variate t distribution} and it will be denoted by
\begin{equation*}
\mathbf{t}=\left(\begin{array}{c}
                    \mathbf{t}_{1} \\
                    \vdots \\
                    \mathbf{t}_{k}
                  \end{array}\right)\sim \mathcal{M}t_{n}(n_{0},\ldots,n_{k};
                  \beta_{1}, \ldots,\beta_{k}).
\end{equation*}

\subsection{Multivector variate Pearson type II}

Recall that if $s^{2}\sim\chi_{\nu}^{2}$ and $\mathbf{x}\sim \mathcal{N}_{n}(\mathbf{0},\mathbf{I})$ are
independent,  then $\mathbf{r}=(s^{2}+||\mathbf{x}||^{2})^{-n/2}\mathbf{x}$ is said to have a
multivariate Pearson type II distribution. Moreover, if $\mathbf{t}$ has a multivariate t distribution,
then $\mathbf{r}=(1+||\mathbf{t}||^{2})^{-n/2}\mathbf{t}$ has also a Person type II distribution, see
\citet{di:67}.

With that remark in mind, assume that $\mathbf{t} \sim \mathcal{M}t_{n}(n_{0},\ldots,n_{k}; \beta_{1},
\ldots,\beta_{k})$ and define $\mathbf{t}_{i}=(1-||\mathbf{r}_{i}||^{2})^{-n/2}\mathbf{r}_{i}$,
$i=1,\ldots,k$, then
$$
  \bigwedge_{i=1}^{k}(d\mathbf{t}_{i}) = \prod_{i=1}^{k}(1-||\mathbf{r}_{i}||^{2})^{-\left(n_{i}/2+1\right)}
   \bigwedge_{i=1}^{k}(d\mathbf{r}_{i}).
$$
And
$$
  dF_{\mathbf{r}_{1},\ldots,\mathbf{r}_{k}}(\mathbf{r}_{1},\ldots,\mathbf{r}_{k}) = \frac{\Gamma\left[n^{*}/2\right]
  \prod_{i=1}^{k}\beta_{i}^{-n_{i}/2}}{\Gamma\left[n_{0}/2\right]\pi^{n^{*}/2}}\hspace{6cm}
$$
$$
  \times \left(1+\sum_{i=1}^{k}\beta_{i}^{-1}||(1-||\mathbf{r}_{i}||^{2})^{-n/2}\mathbf{r}_{i}||^{2}\right)^{-n^{*}/2}
    \prod_{i=1}^{k}(1-||\mathbf{r}_{i}||^{2})^{-\left(n_{i}/2+1\right)}\bigwedge_{i=1}^{k}(d\mathbf{r}_{i}).
$$
Noting that
$$
  1+\sum_{i=1}^{k}\beta_{i}^{-1}||(1-||\mathbf{r}_{i}||^{2})^{-n/2}\mathbf{r}_{i}||^{2}= \prod_{i=1}^{k}
  \left(1-||\mathbf{r}_{i}||^{2}\right)^{-1}\hspace{5cm}
$$
$$
  \hspace{3cm}\times \left(\prod_{i=1}^{k}\left(1-||\mathbf{r}_{i}||^{2}\right)+
  \sum_{i=1}^{k}\prod_{j=1,j\neq i}^{k}\beta_{i}^{-1}\left(1-||\mathbf{r}_{j}||^{2}\right)||\mathbf{r}_{i}||^{2}\right).
$$
Then

$$
   dF_{\mathbf{r}_{1},\ldots,\mathbf{r}_{k}}(\mathbf{r}_{1},\ldots,\mathbf{r}_{k}) =
   \frac{\Gamma\left[n^{*}/2\right) \prod_{i=1}^{k}\beta_{i}^{-n_{i}/2}} {\Gamma\left[n_{0}/2\right]\pi^{n^{*}/2}}
   \prod_{i=1}^{k}(1-||\mathbf{r}_{i}||^{2})^{\sum_{j=1,j\ne i}^{k}\left(n_{i}/2+1\right)-1}
$$
$$
   \times \left(\prod_{i=1}^{k}\left(1-||\mathbf{r}_{i}||^{2}\right)+
  \sum_{i=1}^{k}\prod_{j=1,j\neq i}^{k}\beta_{i}^{-1}\left(1-||\mathbf{r}_{j}||^{2}\right)||\mathbf{r}_{i}||^{2}\right)
    ^{-n^{*}/2} \bigwedge_{i=1}^{k}(d\mathbf{r}_{i}).
$$
This law shall be termed \emph{multivector variate Pearson type II distribution}, with notation
\begin{equation*}
\mathbf{r}=\left(\begin{array}{c}
                    \mathbf{r}_{1} \\
                    \vdots \\
                    \mathbf{r}_{k}
                  \end{array}\right)\sim \mathcal{MP}^{II}_{n}(n_{0},\ldots,n_{k};
                  \beta_{1}, \ldots,\beta_{k}).
\end{equation*}
In a similar way to the multivector variate $\chi^{2}$ generalised t distribution, if we define
$\mathbf{t}_{i}=(1-||\mathbf{r}_{i}||^{2})^{-n/2}\mathbf{r}_{i}$, $i=1,\ldots,k$, then we obtain the so
called \emph{multivector variate $\chi^{2}$ generalised-Pearson type II distribution}:
$$
  dF_{s_{0},\mathbf{r}_{1},\ldots,\mathbf{r}_{k}}(s_{0},\mathbf{r}_{1},\ldots,\mathbf{r}_{k}) =
  \frac{\pi^{n_{0}/2}s_{0}^{n^{*}/2-1}}{\Gamma\left[n_{0}/2\right]\prod_{i=1}^{k}\sigma_{i}^{n_{i}}}
  \hspace{5cm}
$$
$$
  \times h\left(\sigma_{0}^{-2}\left(1+\sigma_{0}^{2}
  \sum_{i=1}^{k}\sigma_{i}^{-2}||(1-||\mathbf{r}_{i}||^{2})^{-n/2}\mathbf{r}_{i}||^{2}\right)s_{0}\right)
$$
$$
   \times \prod_{i=1}^{k}(1-||\mathbf{r}_{i}||^{2})^{-(n_{i}/2+1)}(ds_{0})\wedge
   \left(\bigwedge_{i=1}^{k}(d\mathbf{r}_{i})\right),\hspace{1.5cm}
$$
$$
   = \frac{\pi^{n_{0}/2}s_{0}^{n^{*}/2-1}}{\Gamma\left[n_{0}/2\right]\prod_{i=1}^{k}\sigma_{i}^{n_{i}}}
   h\left(\sigma_{0}^{-2}\left(1+\sigma_{0}^{2} \sum_{i=1}^{k}\sigma_{i}^{-2}(1-||\mathbf{r}_{i}||^{2})
   ||\mathbf{r}_{i}||^{2}\right)s_{0}\right)
$$
$$
  \times\prod_{i=1}^{k}(1-||\mathbf{r}_{i}||^{2})^{-(n_{i}/2+1)}(ds_{0})\wedge\left(\bigwedge_{i=1}^{k}
  (d\mathbf{r}_{i})\right)\hspace{1.5cm}.
$$

Note that the procedure implemented in this section can be generalised and used for deriving an
uncountable number of combinations of families of distributions in the setting of  multivector variate
distributions.

\section{\large Some related distributions}\label{sec:4}

In this section we derive several consequences implicit in the results of the previous apart. In general,
we start with a multivector variate distribution and we are interested in the joint distribution of the
square norms of each subvector, the the multivariate distribution is obtained in such way that the random
components are probabilistically dependent. In the same way a logarithmic distribution and its
generalisation are derived. Moreover, several particular distributions of the literature are derived
straightforwardly.  We also remark that the distributions here derived are of a multivariate type.

\subsection{Multivariate generalised $\chi^{2}$ distribution}\label{subsection:generalizedchi}

Assume that $$\mathbf{x}=(\mathbf{x}_{1}',\ldots,\mathbf{x}_{k}')'\sim
\mathcal{ME}_{n}(n_{1},\ldots,n_{k};
\mathbf{0},\ldots,\mathbf{0};\sigma_{1}^{2}\mathbf{I}_{n_{1}},\ldots,\sigma_{k}^{2}\mathbf{I}_{n_{k}};h),$$
this is
\begin{eqnarray*}
  &&\mathbf{x} =\left(\begin{array}{c}
                    \mathbf{x}_{1} \\
                    \vdots \\
                    \mathbf{x}_{k}
                  \end{array}\right) \in\Re^{n}, \boldsymbol{\mu} = \left(\begin{array}{c}
                    \mathbf{0} \\
                    \vdots \\
                    \mathbf{0}
                  \end{array}\right), \mathbf{x}_{i}\in\Re^{n_{i}},\sum_{i=1}^{k}n_{i}=n,
   \\
   &&\boldsymbol{\Sigma}=\left(
                           \begin{array}{cccc}
                             \sigma_{1}^{2}\mathbf{I}_{n_{1}} & \mathbf{0} & \cdots & \mathbf{0} \\
                             \mathbf{0} & \sigma_{2}^{2}\mathbf{I}_{n_{2}} & \cdots & \mathbf{0} \\
                             \vdots & \vdots & \ddots & \vdots \\
                             \mathbf{0} & \mathbf{0} & \cdots & \sigma_{k}^{2}\mathbf{I}_{n_{k}} \\
                           \end{array}
                         \right),
\end{eqnarray*}
where $\sigma_{i}^{2}>0$ for all $i=1,\ldots,k$.

Then
\begin{equation*}
    dF_{\mathbf{x}_{1},\ldots,\mathbf{x}_{k}}(\mathbf{x}_{1},\ldots,\mathbf{x}_{k})=\prod_{i=1}^{k}(\sigma_{i})^{n_{i}/2}
    h\left(\sum_{i=1}^{k}\sigma_{i}^{-2}||\mathbf{x}_{i}||^{2}\right)\bigwedge_{i=1}^{k}(d\mathbf{x}_{i}).
\end{equation*}
Define $u_{i}=||\mathbf{x}_{i}||^{2}$, then writing $\mathbf{x}_{i}=t_{i}\mathbf{g}_{i},$ $t_{i}>0$,
$i=1,\ldots,k$, $\mathbf{g}_{i}\in \mathcal{V}_{1,n_{i}}$. We have that
$u_{i}=||\mathbf{x}_{i}||^{2}=||t_{i} \mathbf{g}_{i}||^{2}=t_{i}^{2}||\mathbf{g}_{i}||^{2}=t_{i}^{2}$.
Then
\begin{equation*}
    \bigwedge_{i=1}^{k}(d\mathbf{x}_{i})=\bigwedge_{i=1}^{k}\left(2^{-1}u_{i}^{n_{i}/2-1}(du_{i})(\mathbf{g}_{i}'d\mathbf{g}_{i})\right)
    =2^{-k}\prod_{i=1}^{k}u_{i}^{n_{i}/2-1}\bigwedge_{i=1}^{k}(du_{i})\wedge\bigwedge_{i=1}^{k}(\mathbf{g}_{i}'d\mathbf{g}_{i}).
\end{equation*}
Observing that
\begin{eqnarray}
   \int_{\mathbf{g}_{1}\in\mathcal{V}_{1,n_{1}}} \cdots \int_{\mathbf{g}_{k}\in\mathcal{V}_{1,n_{k}}}
    \bigwedge_{i=1}^{k}(\mathbf{g}_{i}'d\mathbf{g}_{i}) &=& \prod_{i=1}^{k}\int_{\mathbf{g}_{i}\in\mathcal{V}_{1,n_{i}}}
    (\mathbf{g}_{i}'d\mathbf{g}_{i})=\prod_{i=1}^{k}\frac{2\pi^{n_{i}/2}}
    {\Gamma\left[n_{i}/2\right]} \nonumber\\
   &=& \frac{2^{k}\pi^{n/2}}{\prod_{i=1}^{k}\Gamma\left[n_{i}/2\right]}, \label{eqint}
\end{eqnarray}
where $n = n_{1}+\cdots+n_{k}$. Then, we finally obtain
\begin{equation}\label{eqmultigenchi}
    dF_{u_{1},\ldots,u_{k}}(u_{1},\ldots,u_{k})=\frac{\pi^{n/2}}
    {\prod_{i=1}^{k}\left(\sigma_{i}^{n_{i}}\Gamma\left[n_{i}/2\right]\right)}
    h\left(\sum_{i=1}^{k}\sigma_{i}^{-2}u_{i}\right)\prod_{i=1}^{k}u_{i}^{n_{i}/2-1}\bigwedge_{i=1}^{k}(du_{i}),
\end{equation}
which shall be termed \emph{multivariate generalised $\chi^{2}$ distribution}, a function denoted by
\begin{equation*}
    \mathbf{u}=\left(
        \begin{array}{c}
          u_{1} \\
          \vdots \\
          u_{k} \\
        \end{array}
      \right)\sim \mathcal{G}\chi_{k}^{2}\left(n_{1},\ldots,n_{k};\sigma_{1}^{2},\ldots,\sigma_{k}^{2};h\right).
\end{equation*}
This sets a proposal for different multivariate distributions, with probabilistic dependent elements, see
\citet{cn:84}, \citet{ln:82}, and \citet{ol:03}. Without any exception, they start with a joint
distribution of independent variables and a set of transformations generate a set of dependent variables.
In our case, we start with a join density function that, except in the Gaussian case, provides a priori
probabilistic dependent variables.

\subsection{Multivariate Beta type I distribution}

Assume that $\mathbf{r}=(\mathbf{r}_{1}',\ldots,\mathbf{r}_{k}')'\sim \mathcal{MP}^{II}_{n}(n_{0},
\ldots,n_{k}; \beta_{1}, \ldots,\beta_{k}),$ and define $b_{i}=||\mathbf{r}_{i}||^{2}$. If we set
$\mathbf{r}_{i}=t_{i}\mathbf{g}_{i},$  $t_{i}>0$, $i=1,\ldots,k$, $\mathbf{g}_{i}\in
\mathcal{V}_{1,n_{i}}$. We obtain
$b_{i}=||\mathbf{r}_{i}||^{2}=||t_{i}\mathbf{g}_{i}||^{2}=t_{i}^{2}||\mathbf{g}_{i}||^{2}=t_{i}^{2}$.
Hence
\begin{equation*}
    \bigwedge_{i=1}^{k}(d\mathbf{r}_{i})=\bigwedge_{i=1}^{k}\left(2^{-1}b_{i}^{n_{i}/2-1}(db_{i})(\mathbf{g}_{i}'d\mathbf{g}_{i})\right)
    =2^{-k}\prod_{i=1}^{k}b_{i}^{n_{i}/2-1}\bigwedge_{i=1}^{k}(db_{i})\wedge\bigwedge_{i=1}^{k}(\mathbf{g}_{i}'d\mathbf{g}_{i}).
\end{equation*}

 Thus

$$
  dF_{b_{1},\ldots,b_{k},\mathbf{g}_{1},\ldots,\mathbf{g}_{k}}(b_{1},\ldots,b_{k},\mathbf{g}_{1},\ldots,\mathbf{g}_{k})
  = \frac{\Gamma\left[n^{*}/2\right]\prod_{i=1}^{k}\beta_{i}^{-n_{i}/2}}{\Gamma\left[n_{0}/2\right]
  \pi^{\sum_{i=1}^{k}n_{i}/2}} \hspace{3cm}
$$
$$
   \hspace{1cm}\times\left[\prod_{i=1}^{k}(1-b_{i})+\sum_{i=1}^{k}\prod_{j=1,j\neq i}^{k}\beta_{i}^{-1}(1-b_{j})b_{i}\right]^{-n^{*}/2}
$$
$$
    \hspace{3cm}\times\prod_{i=1}^{k}(1-b_{i})^{\sum_{j=1,j\neq i}^{k}n_{j}/2-1}
   2^{-k}\prod_{i=1}^{k}b_{i}^{n_{i}/2-1}\bigwedge_{i=1}^{k}(db_{i})\wedge\bigwedge_{i=1}^{k}(\mathbf{g}_{i}'d\mathbf{g}_{i}).
$$

Integrating with respect to $\mathbf{g}_{i}\in \mathcal{V}_{1,n_{i}}$, $i=1,\ldots,k$ by using
(\ref{eqint}), the marginal density of $b_{1},\ldots,b_{k}$ follows
$$
  dF_{b_{1},\ldots,b_{k}}(b_{1},\ldots,b_{k}) = \frac{\prod_{i=1}^{k}\beta_{i}^{-n_{i}/2}}{D_{k}
  \left[n_{0}/2,\ldots,n_{k}/2\right]} \hspace{6cm}
$$
$$
  \hspace{4cm}\times\left[\prod_{i=1}^{k}(1-b_{i})+\sum_{i=1}^{k}\prod_{j=1,j\neq i}^{k}\beta_{i}^{-1}(1-b_{j})b_{i}
   \right]^{-n^{*}/2}
$$
$$
  \hspace{3.5cm} \times\prod_{i=1}^{k}(1-b_{i})^{\sum_{j=1,j\neq i}^{k}n_{j}/2-1}
   \prod_{i=1}^{k}b_{i}^{n_{i}/2-1}\bigwedge_{i=1}^{k}(db_{i}),
$$
where
$$
  D_{k}\left[n_{0}/2,\ldots,n_{k}/2\right]=\frac{\displaystyle\prod_{i=0}^{k}\Gamma\left[n_{i}/2\right]}
  {\Gamma\left[n^{*}/2\right]},
$$
recall that $n^{*} = n_{0}+n_{1}+\cdots+n_{k}$.

The distribution is termed \emph{multivariate beta type I distribution} and it will be denoted as
\begin{equation*}
    \mathbf{b}=\left(
        \begin{array}{c}
          b_{1} \\
          \vdots \\
          b_{k} \\
        \end{array}
      \right)\sim \mathcal{M}\beta_{k}^{I}\left(n_{0},\ldots,n_{k};\beta_{1},\ldots,\beta_{k}\right).
\end{equation*}

\bigskip

This distribution was proposed by \citet{ln:82}. But it was derived from a product of independent Gamma
densities, then  a random vector was constructed with  dependent beta type I random variables expressed
as quotients of the gamma variates. Note that when  $\alpha_{i} = n_{i}/2$ and $\lambda_{i} =
1/\beta_{i}$, $i = 1,\dots , r$ are set in our distribution, the density of \citet{ln:82} is obtained.

\subsection{Multivariate beta type II distribution}

This distribution is also known as the \emph{multivariate F distribution}.

Suppose that $\mathbf{t}=(\mathbf{t}_{1}',\ldots,\mathbf{t}_{k}')'\sim
\mathcal{M}t_{n}(n_{0},\ldots,n_{k}; \beta_{1}, \ldots,\beta_{k})$, this is
\begin{eqnarray*}
    dF_{\mathbf{t}_{1},\ldots,\mathbf{t}_{k}}(\mathbf{t}_{1},\ldots,\mathbf{t}_{k})&=&
    \frac{\Gamma\left[n^{*}/2\right]\prod_{i=1}^{k}\beta_{i}^{-n_{i}/2}}{\Gamma\left[n_{0}/2\right]
  \pi^{\sum_{i=1}^{k}n_{i}/2}}\\&&\times\left(1+\sum_{i=1}^{k}\beta_{i}^{-1}||\mathbf{t}_{i}||^{2}\right)
  ^{-n^{*}/2}\bigwedge_{i=1}^{k}(d\mathbf{t}_{i}).
\end{eqnarray*}
Let $f_{i}=||\mathbf{t}_{i}||^{2}$, $i=1,\ldots,k$. Then by the factorisation
$\mathbf{t}_{i}=w_{i}\mathbf{g}_{i}$, $w_{i}>0$, $\mathbf{g}_{i}\in\mathcal{V}_{1,n_{i}}$, we obtain
$f_{i}=||\mathbf{t}_{i}||^{2}=w_{i}^{2}$. Thus
\begin{equation*}
    \bigwedge_{i=1}^{k}(d\mathbf{t}_{i})=\bigwedge_{i=1}^{k}\left(2^{-1}f_{i}^{n_{i}/2-1}(df_{i})
    (\mathbf{g}_{i}'d\mathbf{g}_{i})\right)=2^{-k}\prod_{i=1}^{k}f_{i}^{n_{i}/2-1}\left(\bigwedge_{i=1}^{k}
    df_{i}\right)\wedge\left(\bigwedge_{i=1}^{k}(\mathbf{g}_{i}'d\mathbf{g}_{i})\right).
\end{equation*}
Using (\ref{eqint}), the marginal distribution of $f_{1},\ldots,f_{k}$ is given by
\begin{eqnarray*}
    dF_{f_{1},\ldots,f_{k}}(f_{1},\ldots,f_{k})&=&
    \frac{\prod_{i=1}^{k}\beta_{i}^{-n_{i}/2}}{D_{k}\left[n_{0}/2,\ldots,n_{k}/2\right]}\\&&\times
    \prod_{i=1}^{k}f_{i}^{n_{i}/2-1}\left(1+\sum_{i=1}^{k}\beta_{i}^{-1}f_{i}\right)^{-n^{*}/2}
    \bigwedge_{i=1}^{k}(df_{i}),
\end{eqnarray*}
a fact which shall be denoted by
$$
  \mathbf{f}=(f_{1},\ldots,f_{k})' \sim \mathcal{M}\beta_{k}^{II}(n_{0}, \ldots,n_{k}, \beta_{1},\ldots,\beta_{k}).
$$

This distribution was also derived by  \citet{ln:82}, starting with independent gamma distributions and
setting dependent quotients which lead to required distribution. They named it as the generalised $F$
distribution. In our context it is obtained by taking $\alpha_{i} = n_{i}/2$ and $\lambda_{i} =
1/\beta_{i}$, $i = 1, \dots,r$.

\subsection{Multivariate generalised gamma-beta type I distribution}

Suppose that $(s_{0},\mathbf{r}_{1}',\ldots,\mathbf{r}_{k}')'$ has a multivector variate generalised
gamma-Pearson type II distribution, then proceeding as in the multivariate beta type I distribution, we
have that
$$
  dF_{s_{0},b_{1}\ldots,b_{k}}(s_{0},b_{1}\ldots,b_{k})= \frac{\pi^{n^{*}/2}s_{0}^{n^{*}
  /2-1}}{\prod_{i=0}^{k} \left[\sigma_{i}^{n_{i}}\Gamma\left[n_{i}/2\right]\right]}
  \prod_{i=1}^{k}(1-b_{i})^{-\left(n_{i}/2+1\right)}b_{i}^{n_{i}/2-1}
$$
$$
  \hspace{3cm}\times h\left[s_{0}\sigma_{0}^{-2}\left(1+\sigma_{0}^{2}\sum_{i=1}^{k}\sigma_{i}^{-2}(1-b_{i})
  b_{i}\right)\right](ds_{0})\wedge\left(\bigwedge_{i=1}^{k}db_{i}\right).
$$
This law will  be denoted as $(s_{0},b_{1},\ldots,b_{k})'\sim
\mathcal{MG}\beta_{k+1}^{I}(n_{0},\ldots,n_{k};\sigma_{0}^{2},\ldots,\sigma_{k}^{2};h)$.

\subsection{Multivariate generalised gamma-beta type II distribution}

Assume that $(s_{0},\mathbf{t}_{1}',\ldots,\mathbf{t}_{k}')'$ has a multivector variate generalised
gamma-Pearson type VII distribution, then if we proceed as in the multivariate beta type II distribution,
then we get that
$$
  dF_{s_{0},f_{1}\ldots,f_{k}}(s_{0},f_{1}\ldots,f_{k})= \frac{\pi^{n^{*}/2}s_{0}^{n^{*}/2-1}}{\prod_{i=0}^{k}
  \left[\sigma_{i}^{n_{i}}\Gamma\left[n_{i}/2\right]\right]}   \prod_{i=1}^{k}f_{i}^{n_{i}/2-1}\hspace{4cm}
$$
$$
   \hspace{5cm}\times h\left[s_{0}\sigma_{0}^{-2}\left(1+\sigma_{0}^{2}\sum_{i=1}^{k}\sigma_{i}^{-2}f_{i}\right)\right]
   (ds_{0})\wedge\left(\bigwedge_{i=1}^{k}df_{i}\right).
$$
A fact that will be denoted by
$$
  (s_{0},f_{1},\ldots,f_{k})'\sim \mathcal{MG}\beta_{k+1}^{II}(n_{0},\ldots,n_{k}; \sigma_{0}^{2},\ldots,\sigma_{k}^{2};h).
$$
This distribution shall be termed \emph{multivariate generalised gamma-beta type II distribution}.

\subsection{Multivariate gamma-log gamma distributions}

Suppose that
$$
 \mathbf{u} \sim  \mathcal{MG}\chi_{k}^{2}(n_{1},\ldots,n_{k_{1}};m_{1},\ldots,m_{k_{2}};\sigma_{1}^{2},\ldots,
  \sigma_{k_{1}}^{2};\delta_{1}^{2},\ldots,\delta_{k_{2}}^{2};h),
$$
where
$$
 \mathbf{u}=
 \left[
 \begin{array}{c}
   \mathbf{u} \\
   \mathbf{w}
 \end{array}
 \right]=
 \left[
 \begin{array}{c}
   u_{1} \\
   \vdots \\
   u_{k_{1}} \\
   w_{1} \\
   \vdots \\
   w_{k_{2}}
 \end{array}
 \right],
$$
and $k_{1}+k_{2}=k$. Then
$$
  dF_{u_{1},\ldots,u_{k_{1}},w_{1},\ldots,w_{k_{2}}}(u_{1},\ldots,u_{k_{1}},w_{1},\ldots,w_{k_{2}})\hspace{6cm}
$$
$$
  \hspace{2cm}=\frac{\pi^{\sum_{i=1}^{k_{1}}n_{i}/2+\sum_{j=1}^{k_{2}}m_{j}/2}}{\prod_{i=1}^{k_{1}}
  \left(\sigma_{i}^{n_{i}}\Gamma\left[n_{i}/2\right]\right)\prod_{i=1}^{k_{2}}\left(\delta_{j}^{m_{j}}
  \Gamma\left[m_{j}/2\right]\right)}\prod_{i=1}^{k_{1}}u_{i}^{n_{i}/2-1}\prod_{i=1}^{k_{2}}w_{j}^{m_{j}/2-1}
$$
$$
   \hspace{4cm}\times h\left(\sum_{i=1}^{k_{1}}\sigma_{i}^{-2}u_{i}+\sum_{j=1}^{k_{2}}\delta_{j}^{-2}w_{j}\right)
   \left(\bigwedge_{i=1}^{k_{1}}du_{i}\right) \wedge\left(\bigwedge_{j=1}^{k_{2}}dw_{j}\right).
$$
Define $y_{j}=\log w_{j}$, $dw_{j}=\exp(y_{j})dy_{j}$, $j=1,\ldots,k_{2}$. Hence
$$
  \left(\bigwedge_{i=1}^{k_{1}}du_{i}\right)
   \wedge\left(\bigwedge_{j=1}^{k_{2}}dw_{j}\right)= \exp\left(\sum_{j=1}^{k_{2}}y_{j}\right)\left(\bigwedge_{i=1}^{k_{1}}du_{i}\right)
   \wedge\left(\bigwedge_{j=1}^{k_{2}}dy_{j}\right).
$$
Then the joint density function of $u_{1},\ldots,u_{k_{1}},y_{1},\ldots,y_{k_{2}}$ is given by
$$
  dF_{u_{1},\ldots,u_{k_{1}},y_{1},\ldots,y_{k_{2}}}(u_{1},\ldots,u_{k_{1}},y_{1},\ldots,y_{k_{2}})\hspace{7cm}
$$
$$
   =\frac{\pi^{\sum_{i=1}^{k_{1}}n_{i}/2+\sum_{j=1}^{k_{2}}m_{j}/2}}{\prod_{i=1}^{k_{1}}\left(\sigma_{i}^{n_{i}}
   \Gamma\left[n_{i}/2\right]\right)\prod_{i=1}^{k_{2}}\left(\delta_{j}^{m_{j}}
   \Gamma\left[m_{j}/2\right]\right)}
   \exp\left(\sum_{j=1}^{k_{2}}\frac{m_{j}y_{j}}{2}\right)\prod_{i=1}^{k_{1}}u_{i}^{n_{i}/2-1}
$$
$$
   \hspace{2cm}\times h\left(\sum_{i=1}^{k_{1}}\sigma_{i}^{-2}u_{i}+\sum_{j=1}^{k_{2}}\delta_{j}^{-2}\exp(y_{i})\right)
   \left(\bigwedge_{i=1}^{k_{1}}du_{i}\right) \wedge\left(\bigwedge_{j=1}^{k_{2}}dy_{j}\right).
$$
This fact will be denoted as
$$
\left(
  \begin{array}{c}
    u_{1} \\
    \vdots \\
    u_{k_{1}} \\
    y_{1} \\
    \vdots \\
    y_{k_{2}} \\
  \end{array}
\right) \sim
\mathcal{MG}\log\mathcal{G}_{k}(n_{1},\ldots,n_{k_{1}};m_{1},\ldots,m_{k_{2}};\sigma_{1}^{2},\ldots,\sigma_{k_{1}}^{2};\delta_{1}^{2},
\ldots,\delta_{k_{2}}^{2};h).
$$
Note that for $k_{2}=0$ we obtain the multivariate gama distribution and for $k_{1}=0$ we have the
multivariate log-gamma distribution.

\section{\large Extended parameter distributions}\label{sec:5}

The reader can check that according to the origin of the distributions, they depends on parameters
$n_{i}$ which literature usually restricts to $n_{i} \in \mathbb{N}$. However, such expressions
(densities) are also valid when the  $n_{i}'s \in \mathbb{N}$ are replaced by the parameters $\alpha_{i}
\in \Re^{+} \quad (\alpha_{i} > 0)$, specifically,  we will replace $n_{i}/2$ by $\alpha_{i}$, $i =
0,\dots,n$ and define $\alpha^{*} = \alpha_{0}+\alpha_{1}+\cdots+\alpha_{k}$. This section defines some
previously derived distributions into the new parametrical space. The implications in the associated
inference can be of interest in some situations, because it enlarges the parametric space avoiding the
complexity of integer optimisations and allowing real parameter estimation.

\subsection{Multivector variate t (or Pearson type VII) distribution}
\begin{small}
\begin{equation}\label{eqmultittgen}
    dF_{\mathbf{t}_{1},\ldots,\mathbf{t}_{k}}(\mathbf{t}_{1},\ldots,\mathbf{t}_{k})=\frac{\Gamma\left[\alpha^{*}\right]
    \prod_{i=1}^{k}\beta_{i}^{-\alpha_{i}}} {\Gamma\left[\alpha_{0}\right]\pi^{\alpha^{*}}}
    \left(1+\sum_{i=1}^{k}\beta_{i}^{-1}||\mathbf{t}_{i}||^{2}\right)^{-\alpha^{*}}
    \bigwedge_{i=1}^{k}(d\mathbf{t}_{i}),
\end{equation}
\end{small}
where $\alpha_{0}>0,\alpha_{i}>0,\beta_{i}>0$, $\mathbf{t}_{i} \in \Re^{n_{i}}$, $i=1,\ldots,k$, with
$n_{0}+\cdots+n_{k} = n$. We refer it as the multivector variate t distribution and it will be denoted by
$$
  \mathbf{t}=(\mathbf{t}_{1}',\ldots,\mathbf{t}_{k}')\sim \mathcal{M}t_{n}(\alpha_{0},\ldots,\alpha_{k}; \beta_{1},
  \ldots,\beta_{k}).
$$

\subsection{Multivector variate Pearson type II distribution}

$$
  dF_{\mathbf{r}_{1},\ldots,\mathbf{r}_{k}}(\mathbf{r}_{1},\ldots,\mathbf{r}_{k}) =
  \frac{\Gamma\left[\alpha^{*}\right] \prod_{i=1}^{k}\beta_{i}^{-\alpha_{i}}} {\Gamma\left[\alpha_{0}\right]
  \pi^{\alpha^{*}}}\prod_{i=1}^{k}\left(1-||\mathbf{r}_{i}||^{2}\right)^{\sum_{j=1,j\ne i}^{k}\left(\alpha_{i}+1\right)-1}
  \hspace{2cm}
$$
$$
 \hspace{1cm}
  \times \left(\prod_{i=1}^{k}\left(1-||\mathbf{r}_{i}||^{2}\right)+ \sum_{i=1}^{k}\prod_{j=1,j\neq i}^{k}
  \beta_{i}^{-1}\left(1-||\mathbf{r}_{j}||^{2}\right)||\mathbf{r}_{i}||^{2}\right)^{-\alpha^{*}}
  \bigwedge_{i=1}^{k}(d\mathbf{r}_{i}),
$$
with $\alpha_{0}>0,\alpha_{i}>0,\beta_{i}>0, i=1,\ldots,k$.
\begin{equation*}
\mathbf{r}=(\mathbf{r}_{1}',\ldots,\mathbf{r}_{k}')\sim \mathcal{MP}^{II}_{n}(\alpha_{0}, \ldots,
\alpha_{k}; \beta_{1}, \ldots,\beta_{k}).
\end{equation*}

\subsection{Multivector variate generalised gamma-Pearson type VII distribution}

If $(s_{0},\mathbf{t}_{1}',\ldots,\mathbf{t}_{k}')'\sim \mathcal{MGP}^{VII}_{n}(\alpha_{0}, \ldots,
\alpha_{k}; \sigma_{0}^{2}, \ldots,\sigma_{k}^{2};h)$, then is density is
$$
  dF_{s_{0},\mathbf{t}_{1},\ldots,\mathbf{t}_{k}}(s_{0},\mathbf{t}_{1},\ldots,\mathbf{t}_{k})=
  \frac{\pi^{\alpha_{0}}}{\Gamma\left[\alpha_{0}\right]\prod_{i=0}^{k}\sigma_{i}^{2\alpha_{i}}}
    s_{0}^{\alpha^{*}-1}\hspace{5cm}
$$
$$
   \hspace{3cm}\times h\left(\sigma_{0}^{-2}\left(1+\sigma_{0}^{2}\sum_{i=1}^{k}\sigma_{i}^{-2}||t_{i}||^{2}
   \right)s_{0}\right)(ds_{0})\wedge \left(\bigwedge_{i=1}^{k}(d\mathbf{t}_{i})\right),
$$
where $\alpha_{i}>0, \sigma_{i}>0, i=1,\ldots,k$.

\subsection{Multivector variate generalised gamma-Pearson type II distribution}

\begin{eqnarray}\label{eqmultigammapearsonIIgen}
  dF_{s_{0},\mathbf{r}_{1},\ldots,\mathbf{r}_{k}}(s_{0},\mathbf{r}_{1},\ldots,\mathbf{r}_{k}) &=&
    \frac{\pi^{\alpha_{0}}s_{0}^{\alpha^{*}-1}}{\Gamma\left[\alpha_{0}\right]\prod_{i=1}^{k}\sigma_{i}^{2\alpha_{i}}}
    \nonumber\\
    &&\times h\left(\sigma_{0}^{-2}\left(1+\sigma_{0}^{2} \sum_{i=1}^{k}\sigma_{i}^{-2}(1-||\mathbf{r}_{i}||^{2})
    ||\mathbf{r}_{i}||^{2}\right)s_{0}\right)\nonumber\\
   &&\times\prod_{i=1}^{k}(1-||\mathbf{r}_{i}||^{2})^{-(\alpha_{i}+1)}(ds_{0})\wedge\left(\bigwedge_{i=1}^{k}(d\mathbf{r}_{i})\right).
\end{eqnarray}
where $\alpha_{i}>0, \sigma_{i}>0, i=1,\ldots,k$. We shall denotes this fact as $(s_{0}, \mathbf{r}_{1}',
\ldots, \mathbf{r}_{k}')'\sim \mathcal{MGP}^{II}_{n}(\alpha_{0},\ldots, \alpha_{k};
\sigma_{0}^{2},\ldots, \sigma_{k}^{2};h)$.

\subsection{Multivariate generalised gamma distribution}\label{subsection:multigamma}

\begin{equation}\label{eqmultigammagen}
    dF_{u_{1},\ldots,u_{k}}(u_{1},\ldots,u_{k})=\frac{\pi^{\sum_{i=1}^{k}\alpha_{i}}}
    {\prod_{i=1}^{k}\left(\sigma_{i}^{2\alpha_{i}}\Gamma\left(\alpha_{i}\right)\right)}
    h\left(\sum_{i=1}^{k}\sigma_{i}^{-2}u_{i}\right)\prod_{i=1}^{k}u_{i}^{\alpha_{i}-1}\bigwedge_{i=1}^{k}(du_{i}),
\end{equation}
with $\alpha_{i}>0, \sigma_{i}>0, i=1,\ldots,k$. This fact shall be denoted as
$\mathbf{u}=(u_{1},\ldots,u_{k})'\sim \mathcal{MG}_{k}\left(\alpha_{1},\ldots, \alpha_{k};
\sigma_{1}^{2}, \ldots,\sigma_{k}^{2};h\right)$.

\subsection{Multivariate beta type I distribution}

For $\mathbf{b}=(b_{1}\ldots,b_{k})'$,
$$
 dF_{b_{1},\ldots,b_{k}}(b_{1},\ldots,b_{k}) = \frac{\prod_{i=1}^{k}\beta_{i}^{-\alpha_{i}}}{D_{k}
 \left[\alpha_{0},\ldots,\alpha_{k}\right]}
 \left[\prod_{i=1}^{k}(1-b_{i})+\sum_{i=1}^{k}\prod_{j=1,j\neq i}^{k}\beta_{i}^{-1}(1-b_{j})b_{i}
 \right]^{-\alpha^{*}}
$$
$$
   \hspace{5cm}\times\prod_{i=1}^{k}(1-b_{i})^{\sum_{j=1,j\neq i}^{k}\alpha_{j}-1}
   \prod_{i=1}^{k}b_{i}^{\alpha_{i}-1}\bigwedge_{i=1}^{k}(db_{i}).
$$
We shall denotes this distribution as, $\mathbf{b}\sim \mathcal{M}\beta_{k}^{I}(\alpha_{0}, \ldots,
\alpha_{k}; \beta_{1},\ldots,\beta_{k})$. In addition, observe that $\alpha_{0}>0$, $\alpha_{i}>0$,
$\beta_{i}>0$, $i=1,\ldots,k$, and
$$
  D_{k} \left[\alpha_{0},\ldots, \alpha_{k}\right] = \frac{\displaystyle \prod_{i=0}^{k} \Gamma
  \left[\alpha_{i}\right]}{\Gamma \left[\alpha^{*}\right]}.
$$

\subsection{Multivariate beta type II distribution}

\begin{eqnarray}\label{eqmultibetaIIgen}
    dF_{f_{1},\ldots,f_{k}}(f_{1},\ldots,f_{k})&=&
    \frac{\prod_{i=1}^{k}\beta_{i}^{-\alpha_{i}}}{D_{k}\left[\alpha_{0},\ldots,\alpha_{k}\right]}
    \prod_{i=1}^{k}f_{i}^{\alpha_{i}-1}\nonumber\\&&\times
    \left(1+\sum_{i=1}^{k}\beta_{i}^{-2}f_{i}\right)^{-\alpha^{*}}
    \bigwedge_{i=1}^{k}(df_{i}).
\end{eqnarray}
where $\alpha_{0}>0, \alpha_{i}>0, \beta_{i}>0, i=1\ldots,k$. And this distribution shall be denoted by
$\mathbf{f}=(f_{1},\ldots,f_{k})'\sim \mathcal{M}\beta_{k}^{II}(\alpha_{0}, \ldots,\alpha_{k};
\beta_{1},\ldots,\beta_{k})$.

\subsection{Multivariate generalised gamma-beta type I distribution}

$$
  dF_{s_{0},b_{1},\ldots,b_{k}}(s_{0},b_{1},\ldots,b_{k})=  \frac{\pi^{\alpha^{*}}s_{0}^{\alpha^{*}-1}}
    {\prod_{i=0}^{k}\left(\sigma_{i}^{2\alpha_{i}}\Gamma[\alpha_{i}]\right)}\prod_{i=1}^{k}(1-b_{i})^{-(\alpha_{i}+1)}
    b_{i}^{\alpha_{i}-1}\hspace{1cm}
$$
$$
  \hspace{2cm}\times h\left[s_{0}\sigma_{0}^{-2}\left(1+\sigma_{0}^{2}\sum_{i=1}^{k}\sigma_{i}^{-2}(1-b_{i})b_{i}\right)\right]
    (ds_{0})\wedge\left( \bigwedge_{i=1}^{k}(db_{i})\right).
$$
with $\alpha_{i}>0, \sigma_{i}>0, i=0\ldots,k$. This distribution shall be denoted as
$(s_{0},b_{1},\ldots,b_{k})'\sim \mathcal{MG}\beta_{k+1}^{I}(\alpha_{0}, \ldots,\alpha_{k};
\sigma_{0}^{2}, \ldots,\sigma_{k}^{2};h)$.

\subsection{Multivariate generalised gamma-beta type II distribution}

If $(s_{0},f_{1},\ldots,f_{k})'\sim \mathcal{MG}\beta_{k+1}^{II}(\alpha_{0}, \ldots,\alpha_{k};
\sigma_{0}^{2}, \ldots,\sigma_{k}^{2};h)$, its density is
$$
  dF_{s_{0},f_{1},\ldots,f_{k}}(s_{0},f_{1},\ldots,f_{k})= \frac{\pi^{\alpha^{*}}s_{0}^{\alpha^{*}-1}}
    {\prod_{i=0}^{k}\left(\sigma_{i}^{2\alpha_{i}}\Gamma[\alpha_{i}]\right)}\prod_{i=1}^{k}f_{i}^{\alpha_{i}-1}
    \hspace{3cm}
$$
$$
  \hspace{3cm}\times h\left[s_{0}\sigma_{0}^{-2}\left(1+\sigma_{0}^{2}\sum_{i=1}^{k}\sigma_{i}^{-2}f_{i}\right)\right]
    (ds_{0})\wedge\left(\bigwedge_{i=1}^{k}(df_{i})\right).
$$
where $\alpha_{i}>0, \sigma_{i}>0, i=0\ldots,k$. .

\subsection{Multivariate generalised gamma-log gamma distribution}
If
$$
  \left[
  \begin{array}{c}
    u_{1} \\
    \vdots \\
    u_{k_{1}} \\
    y_{1} \\
    \vdots \\
    y_{k_{2}}
  \end{array}
  \right] \sim
   \mathcal{MG}\log\mathcal{G}_{k}(\alpha_{1},\ldots,\alpha_{k_{1}};\rho_{1},\ldots,\rho_{k_{1}};\sigma_{1}^{2},\ldots,
   \sigma_{k_{1}^{2}};\delta_{1}^{2},\ldots,\delta_{k_{2}};h),
$$
its density function is given by
$$
  dF_{u_{1},\ldots,u_{k_{1}},y_{1},\ldots,y_{k_{2}}}(u_{1},\ldots,u_{k_{1}},y_{1},\ldots,y_{k_{2}})
  \hspace{6cm}
$$
$$
  \hspace{1cm}= \frac{\pi^{\sum_{i=1}^{k_{1}}\alpha_{i}+\sum_{j=1}^{k_{2}}\rho_{j}}\exp\left(\sum_{j=1}^{k_{2}}\rho_{j}y_{j}\right)}
  {\prod_{i=1}^{k_{1}}\left(\sigma_{i}^{2\alpha_{i}}\Gamma[\alpha_{i}]\right)
  \prod_{j=1}^{k_{2}}\left(\delta_{j}^{2\rho_{j}}\Gamma[\rho_{j}]\right)}\prod_{i=1}^{k_{1}}u_{i}^{\alpha_{i}-1}
$$
$$
   \hspace{4cm}\times h\left(\sum_{i=1}^{k_{1}}\sigma_{i}^{-2}u_{i}+\sum_{j=1}^{k_{2}}\delta_{j}^{-2}e^{y_{j}}\right)
   \left(\bigwedge_{i=1}^{k_{1}}du_{i}\right)\wedge\left(\bigwedge_{j=1}^{k_{2}}dy_{j}\right),
$$
where $\alpha_{i}>0, \sigma_{i}>0, i=1,\ldots,k_{1};\rho_{j}>0,\delta_{j}>0; j=1,\ldots k_{2}$, and
$k_{1}+k_{2}=k$.

Finally, we must quote that the distribution of a subvector is straightforwardly derived from the
corresponding definition in the multivector variate case. In the same way, the marginal of each element
follows from the associated definition in the multivariate case. For example, in the multivector variate
$t$ distribution, all subvector also follows a multivector variate $t$ distribution. In the multivariate
beta case, each subvector has a multivariate beta distribution and each element of the vector has beta
distribution, an so on.

\section{\large Application to dependent financial variables}\label{sec:6}

People and companies are interacting in a market of concern directed by ignorance of the future. The new
techniques of economic growth are based on an speculation in the financial market in which a risk is
assumed in order to obtain certain profitability. Over the years the study of risk has become a
competitive advantage for any operator of the financial market, with this knowledge it is possible to
obtain tools that allow analysing and evaluating those patterns for a more successful prediction of the
future, taking into account the circumstances that may affect the performance of his objectives and
enjoying an advantage over his competitors.

In this section we apply one of the results of the paper to real data in the Colombian financial context.
In this case we consider a dependent probabilistic model between the  Bancolombia preferential stock
(PFBCOLOM)  and the Colombian Security Exchange index (COLCAP Index). Classical studies consider both
variables as independent, but they are clearly highly dependent each other.

The sample consist of $m=61$ pairs of COLCAP Index and PFBCOLOM measured from 02/01/2018 to 02/04/2018.
The random variable for the COLCAP Index will be denoted by $U$, meanwhile $V$ will represent the
PFBCOLOM.

From a financial point of view, there is not a priori marginal or joint distribution for both variables.
The financial usualy studies consider a probabilistic models for each individual variable, then the best
fit for the sample is the starting point for the analysis. After some trials, we found that the
individual fit suggest for $U$ and $V$ a gamma distribution. Recall that the \emph{gamma distribution} is
based on a Gaussian law, and we can define a \emph{generalised gamma distribution} under a general
elliptical model. In particular, we can use an elliptical distribution as the \emph{multivariate Kotz
distribution}, with parameters $r, q$ and $s$, thus, we obtain the \emph{Kotz-gamma distribution}, which
contains the \emph{gamma distribution} when $r=1/2, q=s=1$, see Table \ref{table0}. Then, we are in the
context of Subsection \ref{subsection:generalizedchi}, with a bivariate generalised gamma distribution of
random variables $U$ and $V$.

Note that the addressed multivariate generalised $\chi^{2}$ distribution is valid also when the
$n_{i}/2's$ are positive real numbers. Thus, for enrichment the parametric space and the join fit, we
will set $n_{1}/2=\alpha>0$ and $n_{2}/2=\beta>0$ as real numbers, then we finally will consider the
Subsection \ref{subsection:multigamma} to model our joint dependent data. Both variables  $U$ and $V$ are
defined by the scale and shape real positive parameters, $\sigma_{1}, \alpha$ and $\sigma_{2}, \beta$,
respectively.

Two uses of (\ref{eqmultigammagen}) must be addressed: first, it gives the bivariate density function
which models COLCAP Index and PFBCOLOM; and second, it provides the likelihood function for the parameter
estimation, under a dependence and independence random sample.

Recall that our goal is obtain the MLE (maximum likelihood estimation) of  the interest parameters
assuming that the random variables $V$ and $U$ are probabilistically independent and assuming that
$\mathbf{V} = (v_{1},\ldots,v_{m})'$ and $\mathbf{U} = (u_{1},\ldots,u_{m})'$ are independent random
samples, and; obtain the MLE of  the same parameters, but now assuming that the random variables $V$ and
$U$  are probabilistically dependent and supposing that $\mathbf{V} = (v_{1},\ldots,v_{m})'$ and
$\mathbf{U} = (u_{1},\ldots,u_{m})'$ are dependent random samples.

We start with the maximum likelihood function of $\sigma_{1}, \alpha, \sigma_{2}, \beta$ given the
dependent random samples $\mathbf{V} = (v_{1},\ldots,v_{m})'$, $\mathbf{U} =(u_{1},\ldots,u_{m})'$. In
this case, both variables are measured in the same day, then the sample sizes are equal, i.e. $m=61$. To
apply Table \ref{table0} in the dependent case, we must note that $n=2m(\alpha+\beta)$, then the
likelihood function for a \emph{generalised gamma distributio}n is given by
$$
  \mbox{L}(\sigma_{1}, \alpha, \sigma_{2}, \beta,h|u_{1},\ldots,u_{m},v_{1},\ldots,v_{m}) \hspace{6cm}
$$
$$
  \hspace{1cm}= \frac{\pi^{n/2}\prod_{i=1}^{m}u_{i}^{\alpha-1}\prod_{j=1}^{m}v_{i}^{\beta-1}}
    {\sigma_{1}^{2m\alpha}\sigma_{2}^{2m\beta}\Gamma^{m}\left(\alpha\right)\Gamma^{m}\left(\beta\right)}
    h\left(w\right),
$$
with $w=\sigma_{1}^{-2}\sum_{i=1}^{m}u_{i}+\sigma_{2}^{-2}\sum_{j=1}^{m}v_{i}$. Then, the Log Likelihood
function ($\mathfrak{L}() = \log \mbox{L}()$), for the \emph{Kotz-gamma distribution} with positive
parameters $r, q$ and $s$, is written as
$$
  \mathfrak{L}(\sigma_{1}, \alpha, \sigma_{2}, \beta,r,q,s|u_{1},\ldots,u_{m},v_{1},\ldots,v_{m})
  \hspace{6cm}
$$
$$
  \hspace{0.0cm}= \log(s)+[q+m (\alpha+\beta)-1]\log(r)/s+\log\{\Gamma[m (\alpha+\beta)]\}
$$
$$
  \hspace{1.5cm}-\log\{\Gamma[ (q+m (\alpha+\beta)-1)/s]\}+(\alpha-1) a+(\beta-1) b-m \{2 \alpha \log(\sigma_{1})
$$
$$
  \hspace{0.2cm}+2 \beta \log(\sigma_{2})+\log[\Gamma(\alpha)]+\log[\Gamma(\beta)]\}+(q-1) \log(w)-r w^s,
$$
where $a=\sum_{i=1}^{m}\log(u_{i}), b=\sum_{i=1}^{m}\log(v_{i}), c=\sum_{i=1}^{m}u_{i}$ and $d=\sum_{i=1}^{m}v_{i}$.

Now, the likelihood function under independence is the following product of $m$ \emph{generalised gamma
distributions}:
$$
L(\sigma_{1}, \alpha,h|u_{1},\ldots,u_{m}) =\prod_{i=1}^{m}f_{i}(\sigma_{1},\alpha,u_{i}),
$$
where
$f_{i}(\sigma_{1},\alpha,u_{i})=\pi^{\alpha}\sigma_{1}^{-2\alpha}\Gamma^{-1}(\alpha)u_{i}^{\alpha-1}h(\sigma_{1}^{-2}u_{i})$,
for all $i=1,\ldots,m$. A similar expression is obtained for variable $v_{j}, j=1,\ldots,m$, with
parameters $\sigma_{2}$ and $\beta$.

Then, we find the following particular case of the log likelihood function for the estimation of
$\sigma_{1}$ and $\alpha$ for \emph{Kotz-gamma distribution} with parameters $r,q,s$ and sample
$\mathbf{U} = (u_{1},\ldots,u_{m})'$:
$$
 \mathfrak{L}(\sigma_{1}, \alpha,r,q,s|u_{1},\ldots,u_{m}) \hspace{8cm}
$$
$$
  \hspace{1cm}= m \log(s)+m(q+\alpha-1) \log(r)/s-m \log\{\Gamma[(q+\alpha-1)/s]\}
$$
$$
  \hspace{0.1cm}-2 m (q+\alpha-1) \log(\sigma_{1})+(q+\alpha-2) a-r \sigma_{1}^{-2 s} b
$$
where $a=\sum_{i =1}^{m}\log(u_{i})$ and $b=\sum_{i=1}^{m}u_{i}^{s}$.

A similar expression for $\mathfrak{L}(\sigma_{2}, \beta,r,q,s|v_{1},\ldots,v_{m})$ can be obtained.

Now, finding the MLE of the seven parameters in  the dependent case is a typical problem of optimisation,
which usually have some issues. In particular the initial guess for the starting point is problematic. In
our case, we take advantage that the \emph{multivariate gamma distribution} belongs to \emph{multivariate
Kotz-gamma distribution}, when $r=1/2, q=s=1$, then if we handle suitably such parameters from these
values, then the general likelihood function departures in some controlled sense from the \emph{gamma
distribution}. A second important hint to be point out, considers that the MLE's for the independent and
dependent cases coincides in the Gaussian case. Then we just need to used the well known theory of MLE's
for gamma distribution. For instance, following \citet{cw:84}, the MLE's of the scale parameter $\sigma$
and the shape parameter $\alpha$ of the \emph{gamma distribution} are approximated by:
$$
\hat\alpha=(3-t+\sqrt{(t-3)^{2}+24t})/(12t), \quad \hat\sigma=\sqrt{\sum_{i=1}^{m}u_{i}/(2m\hat\alpha)},
$$
where $t=\log\left(\sum_{i=1}^{m}u_{i}/m\right)-\sum_{i=1}^{m}\log(u_{i})/m$.

In the financial problem under consideration, the MLE's of the COLCAP Index ($U$) are:
$$
\hat\sigma_{1}=0.717296 , \quad \hat\alpha=1472.087461.
$$
And for the PFBCOLOM ($V$), the corresponding MLE's are given by:
$$
\hat\sigma_{2}=3.473827  , \quad \hat\beta=1290.788846.
$$
For a general \emph{multivariate Kotz-gamma distribution}, approximated MLE's are not published, but we
can find them by using the above initial starting point.

In the computations we use several optimisation methods included in the package OPTIMX of R, such as
Nelder-Mead, L-BFGS-B,  nlminb, among others. For the correctness of the above gamma-estimates, we have
applied the referred package in both log likelihood dependent and independent functions indexed by
$r=1/2, q=s=1$, and we find exactly the same estimations.

\begin{figure}[h!]
  \begin{center}
   \includegraphics[width=12cm]{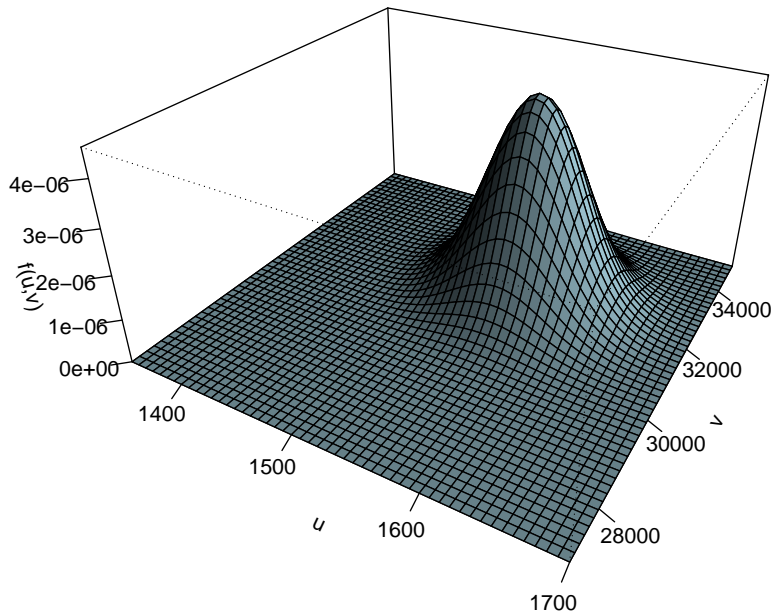}
   \caption{MLE of the Kotz dependent model for the data}\label{fig:dependent}
  \end{center}
\end{figure}

Next we proceed to obtain the MLE  of the seven parameters of the dependent \emph{bivariate Kotz-gamma
distribution} of $u$ and $v$.  Using the gamma-estimates as initial points we have that: $\hat\sigma_{1}=
1.249095$, $\hat\alpha= 1472.083$ $\hat\sigma_{2}= 6.049056$, $\hat\beta= 1290.885$, $\hat{r}=
0.2963556$, $\hat{q}= 1.962216$ and $\hat{s}= 1.129982$, the corresponding \emph{bivariate Kotz-gamma
distribution} is displayed in Figure \ref{fig:dependent}.

\begin{figure}[h!]
  \begin{center}
   \includegraphics[width=12cm]{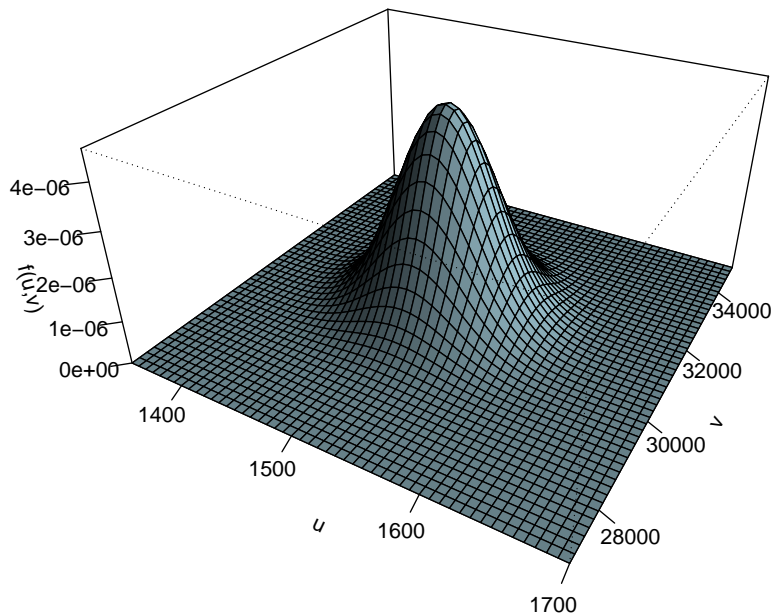}
   \caption{MLE of the Kotz dependent model for the data}\label{fig:independent}
  \end{center}
\end{figure}

Meanwhile, the MLE  of the five parameters of the two independent \emph{univariate Kotz-gamma
distributions} of $u$ and $v$, are given as follows:

For the random variable $u$, we have that $\hat\sigma_{1}= 0.7117504$, $\hat\alpha= 1472.074$,
$\hat{r}=0.5018696$, $\hat{q}= 1.000177$ and $\hat{s}= 0.9978606$. And for the random variable $v$, we
get $\hat\sigma_{2}=  3.482650 $, $\hat\beta=  1290.797$, $\hat{r}=  0.5002108 $, $\hat{q}=  0.9999866$
and $\hat{s}=  1.0005263$. The independent bivariate Kotz gamma distribution is shown in Figure
\ref{fig:independent}.

By visually comparing the values of the estimators of the parameters of interest, it is clear the
importance of establishing whether the random variables under study and their corresponding random
samples are independent or not. It is important to note that it is not interesting to compare
statistically these values differ from the estimators of the parameters under the two contexts, since
they are solutions to different theoretical and/or practical problems.

\section{\large Conclusions}
\begin{enumerate}
  \item The multivector variate distributions, based on multivariate contoured elliptically distribution allow
    to model joint dependent variables, instead of the usual assumption of independence. Moreover, the use of
    elliptical models instead of the classical Gaussian, provides a robust way of modeling  a number of real
    situations. Additionally, note that the distributions obtained, in addition to being used as probabilistic
    models solve the problem of finding the corresponding likelihood functions under independence and dependence.
  \item A real data taking from the Colombian financial context verified the importance to determine if the
    variables under study are independent or not and if the corresponding random samples are also
    independent or not, since the estimators under the two assumptions are not equal.
  \item A number of open problems are of interest, for example, the study of the non central case
    $\boldsymbol{\mu} \neq \mathbf{0}$ (in (\ref{e})) shall allow studied more complex real problems .
\end{enumerate}

\section*{\large Acknowledgements}

This article was partially written during the research stay of the first author, Jos\'e A. D\'{\i}az in
the Department of Agronomy, Division of Life Sciences, Campus Irapuato-Salamanca, University of
Guanajuato, Irapuato, Guanajuato, Mexico.

\end{document}